\documentclass[12pt]{article}
\usepackage{amssymb}
\usepackage{amsfonts}
\usepackage{amsmath}
\usepackage{amsbsy}
\newcommand{\be}{\begin{equation}}
\newcommand{\ee}{\end{equation}}
\newcommand{\bea}{\begin{eqnarray}}
\newcommand{\eea}{\end{eqnarray}}
\newcommand{\ba}{\begin{array}}
\newcommand{\ea}{\end{array}}

\newcommand{\bc}{\begin{center}}
\newcommand{\ec}{\end{center}}
\newcommand{\ben}{\begin{enumerate}}
\newcommand{\een}{\end{enumerate}}
\newcommand{\bfi}{\begin{figure}}
\newcommand{\efi}{\end{figure}}

\newcommand{\bq}{\begin{quote}}
\newcommand{\eq}{\end{quote}}
\newcommand{\bqu}{\begin{quotation}}
\newcommand{\equ}{\end{quotation}}
\newenvironment{emphit}{\begin{itemize}}{\end{itemize}}
\newcommand{\bemp}{\begin{emphit}}
\newcommand{\eemp}{\end{emphit}}

\newcommand{\bt}{\begin{tabular}}
\newcommand{\et}{\end{tabular}}

\newtheorem{myth}{Theorem}[section]

\newtheorem{mycor}{Corollary}[section]

\newtheorem{myrem}{Remark}[section]
\newtheorem{myexam}{Example}[section]

\begin{document}
\date{}
\title{A Double Cryptography Using The Keedwell Cross Inverse Quasigroup
}
\author{T. G. Jaiy\'e\d ol\'a\thanks{All correspondence to be addressed to this author.} \\
Department of Mathematics,\\
Obafemi Awolowo University,\\
Ile Ife 220005, Nigeria.\\
jaiyeolatemitope@yahoo.com\\tjayeola@oauife.edu.ng \and
J. O. Ad\'en\'iran \\
Department of Mathematics,\\
University of Abeokuta, \\
Abeokuta 110101, Nigeria.\\
ekenedilichineke@yahoo.com\\
adeniranoj@unaab.edu.ng} \maketitle
\begin{abstract}
The present study further strenghtens the use of the Keedwell CIPQ
against attack on a system. This is done as follows. The holomorphic
structure of AIPQs(AIPLs) and CIPQs(CIPLs) are investigated.
Necessary and sufficient conditions for the holomorph of a
quasigroup(loop) to be an AIPQ(AIPL) or CIPQ(CIPL) are established.
It is shown that if the holomorph of a quasigroup(loop) is a
AIPQ(AIPL) or CIPQ(CIPL), then the holomorph is isomorphic to the
quasigroup(loop). Hence, the holomorph of a quasigroup(loop) is an
AIPQ(AIPL) or CIPQ(CIPL) if and only if its automorphism group is
trivial and the quasigroup(loop) is a AIPQ(AIPL) or CIPQ(CIPL).
Furthermore, it is discovered that if the holomorph of a
quasigroup(loop) is a CIPQ(CIPL), then the quasigroup(loop) is a
flexible unipotent CIPQ(flexible CIPL of exponent $2$). By
constructing two isotopic quasigroups(loops) $U$ and $V$ such that
their automorphism groups are not trivial, it is shown that $U$ is a
AIPQ or CIPQ(AIPL or CIPL) if and only if  $V$ is a AIPQ or
CIPQ(AIPL or CIPL). Explanations and procedures are given on how
these CIPQs can be used to double encrypt information.
\end{abstract}

\section{Introduction}
\paragraph{}
Let $L$ be a non-empty set. Define a binary operation ($\cdot $) on
$L$ : If $x\cdot y\in L$ for all $x, y\in L$, $(L, \cdot )$ is
called a groupoid. If the system of equations ;
\begin{displaymath}
a\cdot x=b\qquad\textrm{and}\qquad y\cdot a=b
\end{displaymath}
have unique solutions for $x$ and $y$ respectively, then $(L, \cdot
)$ is called a quasigroup. For each $x\in L$, the elements $x^\rho
=xJ_\rho ,x^\lambda =xJ_\lambda\in L$ such that $xx^\rho=e^\rho$ and
$x^\lambda x=e^\lambda$ are called the right, left inverses of $x$
respectively. Now, if there exists a unique element $e\in L$ called
the identity element such that for all $x\in L$, $x\cdot e=e\cdot
x=x$, $(L, \cdot )$ is called a loop.  To every loop $(L,\cdot )$
with automorphism group $AUM(L,\cdot )$, there corresponds another
loop. Let the set $H=(L,\cdot )\times AUM(L,\cdot )$. If we define
'$\circ$' on $H$ such that $(\alpha, x)\circ (\beta,
y)=(\alpha\beta, x\beta\cdot y)$ for all $(\alpha, x),(\beta, y)\in
H$, then $H(L,\cdot )=(H,\circ)$ is a loop as shown in Bruck
\cite{phd82} and is called the Holomorph of $(L,\cdot )$.
\paragraph{}
A loop(quasigroup) is a WIPL(WIPQ) if and only if it obeys the
identity
\begin{equation*}\label{eq:8}
x(yx)^\rho=y^\rho\qquad\textrm{or}\qquad(xy)^\lambda x=y^\lambda.
\end{equation*}
A loop(quasigroup) is a CIPL(CIPQ) if and only if it obeys the
identity
\begin{equation*}\label{eq:8.1}
xy\cdot x^\rho =y\qquad\textrm{or}\qquad x\cdot yx^\rho
=y\qquad\textrm{or}\qquad x^\lambda\cdot
(yx)=y\qquad\textrm{or}\qquad x^\lambda y\cdot x=y.
\end{equation*}
A loop(quasigroup) is an AIPL(AIPQ) if and only if it obeys the
identity
\begin{equation*}
(xy)^\rho=x^\rho y^\rho~or~(xy)^\lambda =x^\lambda y^\lambda
\end{equation*}
Consider $(G,\cdot )$ and $(H,\circ )$ been two distinct
groupoids(quasigroups, loops). Let $A,B$ and $C$ be three distinct
non-equal bijective mappings, that maps $G$ onto $H$. The triple
$\alpha =(A,B,C)$ is called an isotopism of $(G,\cdot )$ onto
$(H,\circ )$ if and only if
\begin{displaymath}
xA\circ yB=(x\cdot y)C~\forall~x,y\in G.
\end{displaymath}
If $(G,\cdot )=(H,\circ )$, then the triple $\alpha =(A,B,C)$ of
bijections on $(G,\cdot )$ is called an autotopism of the
groupoid(quasigroup, loop) $(G,\cdot )$. Such triples form a group
$AUT(G,\cdot )$ called the autotopism group of $(G,\cdot )$.
Furthermore, if $A=B=C$, then $A$ is called an automorphism of the
groupoid(quasigroup, loop) $(G,\cdot )$. Such bijections form a
group $AUM(G,\cdot )$ called the automorphism group of $(G,\cdot )$.

As observed by Osborn \cite{phd89}, a loop is a WIPL and an AIPL if
and only if it is a CIPL. The past efforts of Artzy \cite{phd140,
phd193, phd158, phd30}, Belousov and Tzurkan \cite{phd192} and
present studies of Keedwell \cite{phd176}, Keedwell and Shcherbacov
\cite{phd175, phd177, phd178} are of great significance in the study
of WIPLs, AIPLs, CIPQs and CIPLs, their generalizations(i.e
m-inverse loops and quasigroups, (r,s,t)-inverse quasigroups) and
applications to cryptography.

Interestingly, Adeniran \cite{phd79} and Robinson \cite{phd85},
Oyebo and Adeniran \cite{phd141}, Chiboka and Solarin \cite{phd80},
Bruck \cite{phd82}, Bruck and Paige \cite{phd40}, Robinson
\cite{phd7}, Huthnance \cite{phd44} and Adeniran \cite{phd79} have
respectively studied the holomorphs of Bol loops, central loops,
conjugacy closed loops, inverse property loops, A-loops, extra
loops, weak inverse property loops, Osborn loops and Bruck loops.
Huthnance \cite{phd44} showed that if $(L,\cdot )$ is a loop with
holomorph $(H,\circ)$, $(L,\cdot )$ is a WIPL if and only if
$(H,\circ)$ is a WIPL. The holomorphs of an AIPL and a CIPL are yet
to be studied.

In the quest for the application of CIPQs with long inverse cycles
to cryptography, Keedwell \cite{phd176} constructed the following
CIPQ which we shall specifically call Keedwell CIPQ.
\begin{myth}(Keedwell CIPQ)
Let $(G,\cdot )$ be an abelian group of order $n$ such that $n+1$ is
composite. Define a binary operation '$\circ$' on the elements of
$G$ by the relation $a\circ b=a^rb^s$, where $rs=n+1$. Then
$(G,\circ )$ is a CIPQ and the right crossed inverse of the element
$a$ is $a^u$, where $u=(-r)^3$
\end{myth}
The author also gave examples and detailed explanation and
procedures of the use of this CIPQ for cryptography.

The aim of the present study is to further strenghten the use of the
Keedwell CIPQ against attack on a system. This is done as follows.
\begin{enumerate}
\item The holomorphic structure of AIPQs(AIPLs) and CIPQs(CIPLs) are
investigated. Necessary and sufficient conditions for the holomorph
of a quasigroup(loop) to be an AIPQ(AIPL) or CIPQ(CIPL) are
established. It is shown that if the holomorph of a quasigroup(loop)
is a AIPQ(AIPL) or CIPQ(CIPL), then the holomorph is isomorphic to
the quasigroup(loop). Hence, the holomorph of a quasigroup(loop) is
an AIPQ(AIPL) or CIPQ(CIPL) if and only if its automorphism group is
trivial and the quasigroup(loop) is a AIPQ(AIPL) or CIPQ(CIPL).
Furthermore, it is discovered that if the holomorph of a
quasigroup(loop) is a CIPQ(CIPL), then the quasigroup(loop) is a
flexible unipotent CIPQ(flexible CIPL of exponent $2$).
\item By constructing two isotopic quasigroups(loops) $U$ and $V$ such
that their automorphism groups are not trivial, it is shown that $U$
is a AIPQ or CIPQ(AIPL or CIPL) if and only if  $V$ is a AIPQ or
CIPQ(AIPL or CIPL). Explanations and procedures are given on how
these CIPQs can be used to double encrypt information.
\end{enumerate}

\section{Main Results}
\subsection{Holomorph Of AIPLs And CIPLs}
\begin{myth}\label{3:1}
Let $(L,\cdot )$ be a quasigroup(loop) with holomorph $H(L)$. $H(L)$
is an AIPQ(AIPL) if and only if
\begin{enumerate}
\item $AUM(L)$ is an abelian group,
\item $(\beta^{-1}, \alpha ,I)\in AUT(L)~\forall~\alpha ,\beta\in AUM(L)$ and
\item $L$ is a AIPQ(AIPL).
\end{enumerate}
\end{myth}
{\bf Proof}\\
A quasigroup(loop) is an automorphic inverse property loop(AIPL) if
and only if it obeys the identity
\begin{equation*}
(xy)^\rho=x^\rho y^\rho~or~(xy)^\lambda =x^\lambda y^\lambda.
\end{equation*}
Using either of the definitions of an AIPQ(AIPL) above, it can be
shown that $H(L)$ is a AIPQ(AIPL) if and only if $AUM(L)$ is an
abelian group and $(\beta^{-1}J_\rho, \alpha J_\rho ,J_\rho)\in
AUT(L)~\forall~\alpha ,\beta\in AUM(L)$. $L$ is isomorphic to a
subquasigroup(subloop) of $H(L)$, so $L$ is a AIPQ(AIPL) which
implies $(J_\rho, J_\rho ,J_\rho)\in AUT(L)$. So, $(\beta^{-1},
\alpha ,I)\in AUT(L)~\forall~\alpha ,\beta\in AUM(L)$.

\begin{mycor}\label{3:2}
Let $(L,\cdot )$ be a quasigroup(loop) with holomorph $H(L)$. $H(L)$
is a CIPQ(CIPL) if and only if
\begin{enumerate}
\item $AUM(L)$ is an abelian group,
\item $(\beta^{-1}, \alpha ,I)\in AUT(L)~\forall~\alpha ,\beta\in AUM(L)$ and
\item $L$ is a CIPQ(CIPL).
\end{enumerate}
\end{mycor}
{\bf Proof}\\
A quasigroup(loop) is a CIPQ(CIPL) if and only if it is a WIPQ(WIPL)
and an AIPQ(AIPL). $L$ is a WIPQ(WIPL) if and only if $H(L)$ is a
WIPQ(WIPL).

If $H(L)$ is a CIPQ(CIPL), then $H(L)$ is both a WIPQ(WIPL) and a
AIPQ(AIPL) which implies 1., 2., and 3. of Theorem~\ref{3:1}. Hence,
$L$ is a CIPQ(CIPL). The converse follows by just doing the reverse.

\begin{mycor}\label{3:3}
Let $(L,\cdot )$ be a quasigroup(loop) with holomorph $H(L)$. If
$H(L)$ is an AIPQ(AIPL) or CIPQ(CIPL), then $H(L)\cong L$.
\end{mycor}
{\bf Proof}\\
By 2. of Theorem~\ref{3:1}, $(\beta^{-1}, \alpha ,I)\in
AUT(L)~\forall~\alpha ,\beta\in AUM(L)$ implies $x\beta^{-1}\cdot
y\alpha =x\cdot y$ which means $\alpha =\beta =I$ by substituting
$x=e$ and $y=e$. Thus, $AUM(L)=\{I\}$ and so $H(L)\cong L$.

\begin{myth}\label{3:3.1}
The holomorph of a quasigroup(loop) $L$ is a AIPQ(AIPL) or
CIPQ(CIPL) if and only if $AUM(L)=\{I\}$ and $L$ is a AIPQ(AIPL) or
CIPQ(CIPL).
\end{myth}
{\bf Proof}\\
This is established using Theorem~\ref{3:1}, Corollary~\ref{3:2} and
Corollary~\ref{3:2}.

\begin{myth}\label{3:4}
Let $(L,\cdot )$ be a quasigroups(loop) with holomorph $H(L)$.
$H(L)$ is a CIPQ(CIPL) if and only if $AUM(L)$ is an abelian group
and any of the following is true for all $x,y\in L$ and $\alpha
,\beta\in AUM(L)$:
\begin{enumerate}
\item $(x\beta\cdot y)x^\rho=y\alpha$.
\item $x\beta\cdot yx^\rho=y\alpha$.
\item $(x^\lambda\alpha^{-1}\beta\alpha\cdot y\alpha )\cdot x=y$.
\item $x^\lambda\alpha^{-1}\beta\alpha\cdot (y\alpha \cdot x)=y$.
\end{enumerate}
\end{myth}
{\bf Proof}\\
This is achieved by simply using the four equivalent identities that
define a CIPQ(CIPL):
\begin{displaymath}
xy\cdot x^\rho =y\qquad\textrm{or}\qquad x\cdot yx^\rho
=y\qquad\textrm{or}\qquad x^\lambda\cdot
(yx)=y\qquad\textrm{or}\qquad x^\lambda y\cdot x=y.
\end{displaymath}

\begin{mycor}\label{3:5}
Let $(L,\cdot )$ be a quasigroups(loop) with holomorph $H(L)$. If
$H(L)$ is a CIPQ(CIPL) then the following are equivalent to each
other
\begin{enumerate}
\item $(\beta^{-1}J_\rho, \alpha J_\rho ,J_\rho)\in AUT(L)~\forall~\alpha
,\beta\in AUM(L)$.
\item $(\beta^{-1}J_\lambda, \alpha J_\lambda ,J_\lambda )\in AUT(L)~\forall~\alpha
,\beta\in AUM(L)$.
\item $(x\beta\cdot y)x^\rho=y\alpha$.
\item $x\beta\cdot yx^\rho=y\alpha$.
\item $(x^\lambda\alpha^{-1}\beta\alpha\cdot y\alpha )\cdot x=y$.
\item $x^\lambda\alpha^{-1}\beta\alpha\cdot (y\alpha \cdot x)=y$.
\end{enumerate}
Hence,
\begin{displaymath}
(\beta, \alpha ,I),(\alpha, \beta ,I),(\beta,I, \alpha ),
(I,\alpha,\beta)\in AUT(L)~\forall~\alpha ,\beta\in AUM(L).
\end{displaymath}
\end{mycor}
{\bf Proof}\\
The equivalence of the six conditions follows from Theorem~\ref{3:4}
and the proof of Theorem~\ref{3:1}. The last part is simply.

\begin{mycor}\label{3:6}
Let $(L,\cdot )$ be a quasigroup(loop) with holomorph $H(L)$. If
$H(L)$ is a CIPQ(CIPL) then, $L$ is a flexible unipotent
CIPQ(flexible CIPL of exponent $2$).
\end{mycor}
{\bf Proof}\\
It is observed that $J_\rho =J_\lambda =I$. Hence, the conclusion
follows.

\begin{myexam}\label{3:7}
Let $(L,\cdot )$ be an abelian group with $\textrm{Inn}_\rho
(L)$-holomorph $H(L)$. $H(L)$ is an abelian group.
\end{myexam}
{\bf Proof}\\
In an extra loop $L$, $\textrm{Inn}_\rho (L)=\textrm{Inn}_\lambda
(L)\le AUM(L)$ is a boolean group, hence it is abeilan group. An
abelian group is a commutative extra loop. A commutative extra loop
is a CIPL. So by Corollary~\ref{3:2}, $H(L)$ is a CIPL. $H(L)$ is a
group since $L$ is a group. A group is a CIPL if and only it is
abelian. Thus, $H(L)$ is an abelian group.

\begin{myrem}
The holomorphic structure of loops such as extra loop, Bol-loop,
C-loop, CC-loop and A-loop have been found to be characterized by
some special types of automorphisms such as
\begin{enumerate}
\item Nuclear automorphism(in the case of Bol-,CC- and extra loops),
\item central automorphism(in the case of central and A-loops).
\end{enumerate}
By Theorem~\ref{3:1} and Corollary~\ref{3:2}, the holomorphic
structure of AIPLs and CIPLs is characterized by commutative
automorphisms. The abelian group in Example~\ref{3:7} is a boolean
group.
\end{myrem}

\subsection{A Pair Of AIPLs And CIPLs}
\begin{myth}\label{1:4}
Let $U=(L,\oplus)$ and $V=(L,\otimes )$ be quasigroups such that
$AUM(U)$ and $AUM(V)$ are conjugates in $SYM(L)$ i.e there exists a
$\psi\in SYM(L)$ such that for any $\gamma\in AUM(V)$, $\gamma
=\psi^{-1}\alpha\psi$ where $\alpha\in AUM(U)$. Then, $H(U)\cong
H(V)$ if and only if $x\delta\otimes y\gamma =(x\beta\oplus
y)\delta~\forall~x,y\in L,~\beta\in AUM(U)$ and some
$\delta,\gamma\in AUM(V)$. Hence:
\begin{enumerate}
\item $\gamma\in AUM(U)$ if and only if $(I,\gamma ,\delta )\in
AUT(V)$.
\item if $U$ is a loop, then;
\begin{enumerate}
\item ${\cal L}_{e\delta}\in AUM(V)$.
\item $\beta\in AUM(V)$ if and only if ${\cal R}_{e\gamma}\in
AUM(V)$.
\end{enumerate}
where $e$ is the identity element in $U$ and ${\cal L}_x$, ${\cal
R}_x$ are respectively the left and right translations mappings of
$x\in V$.
\item if $\delta =I$, then $|AUM(U)|=|AUM(V)|=3$ and so
$AUM(U)$ and $AUM(V)$ are boolean groups.
\item if $\gamma =I$, then $|AUM(U)|=|AUM(V)|=1$.
\end{enumerate}
\end{myth}
{\bf Proof}\\
\begin{enumerate}
\item Let $H(L,\oplus )=(H,\circ )$ and $H(L,\otimes )=(H,\odot )$.
$H(U)\cong H(V)$ if and only if there exists a bijection $\phi
~:~H(U)\to H(V)$ such that $[(\alpha ,x)\circ (\beta ,y)]\phi
=(\alpha,x)\phi\odot (\beta ,y)\phi$. Define $(\alpha,x)\phi
=(\psi^{-1}\alpha\psi ,x\psi^{-1}\alpha\psi )~\forall~(\alpha ,x)\in
(H,\circ )$ where $\psi\in SYM(L)$.
\item $H(U)\cong H(V)\Leftrightarrow (\alpha\beta ,x\beta\oplus y)\phi
=(\psi^{-1}\alpha\psi ,x\psi^{-1}\alpha\psi )\odot
(\psi^{-1}\beta\psi ,y\psi^{-1}\beta\psi )\Leftrightarrow
(\psi^{-1}\alpha\beta\psi ,(x\beta\oplus y)\psi^{-1}\alpha\beta\psi
)=(\psi^{-1}\alpha\beta\psi ,x\psi^{-1}\alpha\beta\psi\otimes
y\psi^{-1}\beta\psi )\Leftrightarrow (x\beta\oplus
y)\psi^{-1}\alpha\beta\psi=x\psi^{-1}\alpha\beta\psi\otimes
y\psi^{-1}\beta\psi\Leftrightarrow x\delta\otimes y\gamma
=(x\beta\oplus y)\delta$ where $\delta =\psi^{-1}\alpha\beta\psi$,
$\gamma =\psi^{-1}\beta\psi$.
\item Note that, $\gamma{\cal L}_{x\delta}=L_{x\beta}\delta$ and
$\delta{\cal R}_{y\gamma}=\beta R_y\delta~\forall~x,y\in L$. So,
when $U$is a loop, $\gamma{\cal L}_{e\delta}=\delta$ and
$\delta{\cal R}_{e\gamma}=\beta\delta$. These can easily be used to
prove the remaining part of the theorem.
\end{enumerate}

\begin{myth}\label{1:6}
Let $U=(L,\oplus)$ and $V=(L,\otimes )$ be quasigroups(loops) that
are isotopic under the triple of the form $(\delta^{-1}\beta
,\gamma^{-1},\delta^{-1})$ for all $\beta\in AUM(U)$ and some
$\delta,\gamma\in AUM(V)$ such that their automorphism groups are
non-trivial and are conjugates in $SYM(L)$ i.e there exists a
$\psi\in SYM(L)$ such that for any $\gamma\in AUM(V)$, $\gamma
=\psi^{-1}\alpha\psi$ where $\alpha\in AUM(U)$. Then, $U$ is a AIPQ
or CIPQ(AIPL or CIPL) if and only if  $V$ is a AIPQ or CIPQ(AIPL or
CIPL).
\end{myth}
{\bf Proof}\\
Let $U$ be an AIPQ or CIPQ(AIPL or CIPL), then since $H(U)$ has a
subquasigroup(subloop) that is isomorphic to $U$ and that
subquasigroup(subloop) is isomorphic to a subquasigroup(subloop) of
$H(V)$ which is isomorphic to $V$, $V$ is a AIPQ or CIPQ(AIPL or
CIPL). The proof for the converse is similar.

\subsection{Application To Cryptography}
Let the Keedwell CIPQ be the quasigroup $U$ in Theorem~\ref{1:4}.
Definitely, its automorphism group is non-trivial because as shown
in Theorem~2.1 of Keedwell \cite{phd176}, for any CIPQ, the mapping
$J_\rho~:~x\to x^\rho$ is an automorphism. This mapping will be
trivial only if $U$ is unipotent. For instance, in  Example 2.1 of
Keedwell \cite{phd176}, the CIPQ $(G,\circ )$ obtained is unipotent
because it was constructed using the cyclic group $C_5=<c:~c^5=e>$
and defined as $a\circ b=a^3b^2$. But in Example~2.2, the CIPQ
gotten is not unipotent as a result of using the cyclic group
$C_{11}=<c:~c^{11}=e>$. Thus the choice of a Keedwell CIPQ which
suits our purpose in this work for a cyclic group of order $n$ is
one in which $rs=n+1$ and $r+s\ne n$. Now that we have seen a sample
for the choice of $U$, the quasigroup $V$ can then be obtained as
shown in Theorem~\ref{1:4}. By Theorem~\ref{1:6}, $V$ is a CIPQ.

In Keedwell \cite{phd176}, the author's method of application is as
follows. It is assumed that the message to be transmitted can be
represented as single element $x$ of the quasigroup $U$ and that
this is enciphered by multiplying by another element $y$ of $U$ so
that the encoded message is $yx$. At the receiving end, the message
is deciphered by multiplying by the inverse of $y$. Now, according
to Theorem~\ref{1:4}, by the choice of the mappings $\alpha
,\beta\in AUM(U)$ and $\psi\in Sym(L)$ to get the mappings
$\delta,\gamma$, a CIPQ $V$ can be produced following
Theorem~\ref{1:4}. So, the secret keys for the systems are $\{\alpha
,\beta,\psi\}\equiv\{\delta,\gamma\}$. Thus whenever a set of
infomation or messages is to be transmitted, the sender will
enciphere in the Keedwell CIPQ(as described earlier on) and then
enciphere again with $\{\alpha ,\beta,\psi\}\equiv\{\delta,\gamma\}$
to get a CIPQ $V$ which is the set of encoded messages. At the
receiving end, the message $V$ is deciphered by using an inverse
isotopism(i.e inverse key $\{\alpha
,\beta,\psi\}\equiv\{\delta,\gamma\}$) to get $U$ and then deciphere
again(as described earlier on) to get the messages. The secret key
can be changed over time. The method described above is a double
encryption and its a double protection. It protects each piece of
information(element of the quasigroup) and protects the combined
information(the quasigroup as a whole). Its like putting on a pair
of socks and shoes or putting on under wears and clothes, the body
gets better protection.

\end{document}